\documentclass[11pt]{article}

\def\N{\mathbb{ N}}

\def\Z{\mathbb{Z}}

\def\Par{\mathcal{P}}

\usepackage{epsfig}
\usepackage{amssymb,amsfonts,latexsym}
\newtheorem{theo}{Theorem}

\newtheorem{coro}{Corollary}

\newtheorem{lem}{Lemma}
\newenvironment{proof}{{\em Proof.}}{\mbox{}\hfill $\Box$\medskip}
\newenvironment{ex}{\noindent{\bf Example. }}{\smallskip}
\newenvironment{rem}{\noindent{\bf Remark. }}{\smallskip}
\newcommand{\la}{\lambda}
\title{Some new identities for Schur functions }
\author{Fr\'ed\'eric Jouhet and  Jiang Zeng\\
\small Institut Girard Desargues, Universit\'e Claude Bernard
(Lyon 1)\\ 
\small 43, bd du 11 Novembre 1918, 69622 Villeurbanne
Cedex, France\\ 
\small \texttt{jouhet@desargues.univ-lyon1.fr,
zeng@desargues.univ-lyon1.fr}\\
\date{Dedicated to Dominique Foata on the occasion of his 65th birthday}\\
}
\begin{document}
\maketitle
\begin{abstract}
Some new identities for  Schur functions are proved. In particular, 
we settle in the affirmative 
a recent conjecture of Ishikawa-Wakayama~\cite{IW} 
and solve a problem
raised by Bressoud~\cite{Br2}.
\end{abstract}

\section{Introduction}
We fix a positive integer $n$ and let $X=(x_1,\ldots ,x_n)$ be a set of $n$ 
independent variables.
For each partition 
$\la=(\la_1\geq \la_2\geq \cdots \geq \la_n\geq 0)$
of length $\leq n$, the 
Schur function  $s_\la(X)$ are usually defined  as follows~\cite{Ma}: 
$$
s_\la(X)=\det\left(x_i^{\la_j+n-j}\right)_{1\leq i,j\leq
n}/\det\left(x_i^{n-j}\right)_{1\leq i,j\leq n}. 
$$
In this paper we shall follow the standard definitions and
notations of Macdonald's book~\cite{Ma}.  Thus the 
\emph{Ferrers diagram} of $\la$ is the subset 
 $\{(i,\,j)|j\geq 1, i\leq \la_j\}$ of $\N^2$.
If the diagram of $\mu$ is included in
that of $\la$ we note
$\mu\subseteq \la$ and the skew diagram $\la/\mu$ is called  a
\emph{horizontal strip} (or h.s. for short) if there is at most one cell
in each column of  $\la/\mu$. For any partition $\lambda$ we note
$c_j:=c_j(\la)$ the number of columns of length $j$ in $\la$, i.e.
$c_j=\la_j-\la_{j+1}$ and define 
$$
f_\la(a,b)=\prod_{j\,
\hbox{\scriptsize odd}}\frac{a^{c_j+1}-b^{c_j+1}}{a-b}\,
\prod_{j\, \hbox{\scriptsize even}}\frac{1-(ab)^{c_j+1}}{1-ab}. 
$$
Since $f_\la(a,0)=a^{c(\la)}$, where $c(\la)$ is the number of 
columns of odd length of $\la$, a  
classical identity of Littlewood~\cite{Ma} reads then as follows :
\begin{equation}\label{litt}
\sum_{\la}f_\la(a,0)s_\la(X)
=\prod_{i}(1-ax_i)^{-1}\prod_{j<k}(1-x_jx_k
)^{-1}.
\end{equation}
Set
$$
\Phi(X;a,b):=\prod_{i}(1-ax_i)^{-1}
(1-bx_i)^{-1}\prod_{j<k}(1-x_jx_k)
^{-1}.
$$ 
In a recent paper~\cite{IW}, Ishikawa and Wakayama  gave the following
extension of (1):
\begin{equation} \label{iw}
\sum_{\la}f_\la(a,b)\,s_\la(X) =\Phi(X;a,b).
\end{equation}
As pointed  by Bressoud~\cite{Br2}, when $(a,b)=(1,0)$, $(1,-1)$ and
$(0,0)$, identity (\ref{iw}) reduces to the following interesting
known identities respectively:
\begin{eqnarray}
\sum_{\la} s_\la(X)&=&\prod_{i=1}^n{1\over 1-x_i}\prod_{1\leq
i<j\leq n}{1\over 1-x_ix_j},\\ \sum_{\la\,\textrm{even}}
s_\la(X)&=&\prod_{i=1}^n{1\over 1-x_i^2} \prod_{1\leq i<j\leq
n}{1\over 1-x_ix_j},\\ \sum_{\la'\,\textrm{even}}
s_\la(X)&=&\prod_{1\leq i<j\leq n}{1\over 1-x_ix_j},
\end{eqnarray}
where $\la'$ is the conjugate  of $\la$.

In this paper we shall give two generalizations of 
Ishikawa and Wakayama's formula (\ref{iw}). To state them we need some definition.

For $r\geq 0$,
let  $h_r(X)$ (resp. $e_r(X)$) be the \emph{homogeneous} 
(resp. \emph{elementary}) symmetric function of
$X$ and set
\begin{eqnarray*}
P_r(a,b,c)&=&\sum_{k=0}^r\frac{a^{k+1}-b^{k+1}}{a-b}\,
\frac{1-(ab)^{r-k+1}}{1-ab}c^k,\\
Q_r(a,b,c)&=&\sum_{k=0}^rh_{r-k}(a,b,c)(abc)^k.
\end{eqnarray*}
For any positive integer sequence $\xi=(\xi_1,\ldots, \xi_r,\ldots)$, where
$\xi_r\neq 0$ for only a finite number of integers $r$, set
$$ 
F_{\xi}(a,b,c)=h_{\xi_1}(a,b,c)\prod_{k\geq
1}P_{\xi_{2k}}(a,b,c)Q_{\xi_{2k+1}}(a,b,c). 
$$ 
For any integer $i\geq 1$, let
 $\varepsilon_i$ be the $i^{th}$ vector of the canonical
basis of $\Z^\infty$ and 
introduce  
the operator $\delta_i: \delta_i\xi=\xi-\varepsilon_i-\varepsilon_{i+1}$
for  $\xi\in  \N^\infty$. 
Set  $\delta_i F_{\xi}(a,b,c)=F_{\delta_i\xi}(a,b,c)$, where 
$P_k=Q_k=0$ if $k<0$ by convention.
 Hence, to  any partition $\la$ of length $\leq n$ we can associate 
the polynomial
$$
f_\la(a,b,c):=\sum_{k=0}^n(-abc)^k\sum_{i_1<\cdots <i_k}
\delta_{i_1}\cdots \delta_{i_k}
F_{\Gamma(\la)}(a,b,c),
$$ 
where $\Gamma(\la)=(c_1,c_2,\ldots)$ is 
the sequence of the multiplicities of parts in the dual of  $\la$,
or $c_j$ is the number of columns of length $j$ in $\la$.

Now we can state our first generalization of (\ref{iw}), which gives in fact
a positive answer
 to a  conjecture of Ishikawa and Wakayama~{\cite{IW}.
\begin{theo} We have
$$ 
\sum_\la
f_\la(a,b,c)s_\la(X)=\Phi(X;a,b)\prod_{i}(1-cx_i)^{-1}. 
$$
\end{theo}

On the other hand,  
Macdonald~\cite[ p. 83-84]{Ma}, 
D\'esarm\'enien-Stembridge~\cite{De,St} and Okada~\cite{Ok} have
given 
\emph{bounded versions} of identities (3)-(5), respectively, as follows~:

\begin{theo}[Macdonald] For 
 non negative integers $m$ and $n$, 
$$
\sum_{\la_1\leq m}s_\la(X)=
{\det\left(x_i^{j-1}-x_i^{m+2n-j}\right)\over
\prod_{i=1}^n(1-x_i)\prod_{i<j}(x_i-x_j)(x_ix_j-1)}.
$$
\end{theo}

\begin{theo}[D\'esarm\'enien-Stembridge]\label{th:DS} For 
 non negative integers $m$ and $n$,
$$ \sum_{\la_1\leq 2 m\atop \la \,\hbox{\scriptsize
even}}s_\la(X) ={\det\left(x_i^{j-1}-x_i^{2m+2n+1-j}\right)\over
\prod_{i=1}^n(1-x_i^2) \prod_{i<j}(x_ix_j-1)(x_i-x_j)}. $$
\end{theo}

\begin{rem}
This result follows immediately from Macdonald's formula.
Indeed Pieri's formula implies:
 $$
\sum_{k=0}^ne_k(X)\, \sum_{\la_1\leq 2 m\atop \la
\hbox{\scriptsize even}}s_\la(X)= \sum_{\la_1\leq 2 m+1}s_\la(X).
$$ 
Since $\sum_{k=0}^ne_k(X)=\prod_{i=1}^n(1+x_i)$, we get
immediately theorem~\ref{th:DS} by applying Macdonald's formula.
\end{rem}

\begin{theo}[Okada] For non negative  integers $m$ and $n$ which is even,
$$ 
\sum_{\la_1\leq m\atop \la'\,
\hbox{\scriptsize even}}s_\la(X) =
{1\over 2}{\det\left(x_i^{j-1}-x_i^{m+2n-1-j}\right)+
\det\left(x_i^{j-1}+x_i^{m+2n-1-j}\right)\over
\prod_{i<j}(x_ix_j-1)(x_i-x_j)}.
 $$
\end{theo}

After giving elementary proofs of (2) and of the last three
identities \cite{Br1,Br2},  Bressoud~\cite{Br2} raised the problem 
of finding an extension of (2)  for bounded partitions. 
Our  second generalization of 
(2) will give an answer to Bressoud's problem~\cite{Br2}.

For any sequence $\xi\in \{\pm 1\}^n$, we
denote by $|\xi|_{-1}$ the number of $-1$'s in the sequence $\xi$, set
$X^\xi=\{x_1^{\xi_1},\ldots, x_n^{\xi_n}\}$ and 
$$
D(\xi,z)=1-z\prod_ix_i^{(\xi_i-1)/2}.
$$

\begin{theo} For non negative  integers $m$ and $n$,
$$
\sum_{\la\subseteq  (m^n)}f_\la(a,b)s_\la(X)
=\sum_{\xi\in \{\pm 1\}^n}\beta(\xi,a,b)
\Phi(X^{\xi}; a,b)\prod_i x_i^{m(1-\xi_i)/2}
$$
where the coefficient $\beta(\xi,a,b)$ is equal to
$$ 
\left\lbrace\begin{array}{ll}\displaystyle 
\left({a^{m+1}\over D(\xi,1/a)}
-{b^{m+1}\over  D(\xi,1/b) }\right)
{D(\xi,a)D(\xi,b)\over a-b}& \textrm{if} \; |\xi|_{-1} \; \textrm{odd} ,\\
&\\
\displaystyle \left({1\over D(\xi,1)}
-{(ab)^{m+1}\over  D(\xi,1/ab) }\right)
{D(\xi,1)D(\xi,ab)\over 1-ab}& \textrm{if}\;  |\xi|_{-1}\; \textrm{even} .
		       \end{array}
\right.
$$
\end{theo}
\begin{rem}
Assume that $|a|<1$, $|b|<1$ and $|x_i|<1$ ($1\leq i\leq n$) and let
$m\to \infty$, then all the summands tend to 0 except the one
corresponding to $|\xi|_{-1}=0$, which tends to $\Phi(X;a,b)$.
Therefore theorem~5 reduces to (2) when $m\to \infty$. 
\end{rem}

We shall give the proof of 
theorem~1 in section~2 and that of theorem~5 in section~3 using Macdonald's 
approach~\cite{Ma}.
Finally, in section~4, we will show that when $(a,b)=(1,0)$, $(1,-1)$ and
$(0,0)$, theorem~5 reduces actually to the above results 
of Macdonald, D\'esarm\'enien-Stembridge and Okada respectively.

\section{Proof of theorem~1}
Let $\Par$ be the set of partitions of length $\leq n$.
Given a partition $\la\in \Par$, we note
$H(\la)$  the set of partitions $\mu\in \Par$ 
such that $\la/\mu$ is a horizontal strip.
As noticed at the end of \cite{IW}, identity (\ref{iw}) can be derived from
Littlewood's formula
(\ref{litt}) and the so-called Pieri formula (see \cite{Ma}):
\begin{equation}
s_\mu(X)\,h_k(X)=\sum_{\la:\mu\in H(\la) \atop |\la/\mu|=k}s_\la(X). 
\label{eq:pieri}
\end{equation}
In the same vain, we shall 
derive  theorem~1 from (\ref{iw}) and (\ref{eq:pieri}).
We first review such a proof for  (\ref{iw}).
By virtue of  (1) and (\ref{eq:pieri}), identity  (2) is 
equivalent to the following:
\begin{equation} \label{iw'}
f_\la(a,b)=\sum_{\mu\in H(\la)}b^{|\la/\mu|} a^{c(\mu)}.
\end{equation}
Let  $A_j(\la)$ be the 
subdiagram of $\la$ consisting of
 $c_j$ columns of length $j$ for $j\geq 1$.
Thus choosing a partition $\mu$ in $H(\la)$ is equivalent to
choose $r$ left-most (resp. the rest $c_j-r$ ) 
columns of length $j$ (resp. $j-1$ ) for $\mu$ within each block $A_j(\mu)$.
Clearly  the corresponding weight is 
$$\left\lbrace\begin{array}{ll}
\sum_{r=0}^{c_j}a^{c_j-r}b^{r}=\frac{a^{c_j+1}-b^{c_j+1}}{a-b}\quad
&\textrm{if  $j$ is odd}, \\
&\\
\sum_{r=0}^{c_j}(ab)^{r}=\frac{1-(ab)^{c_j+1}}{1-ab}
&\textrm{if $j$ is even}.
	      \end{array}
\right.
$$
Multiplying the weights on all $j\geq 1$ yields (\ref{iw'}).

Each pair $(\la,\mu)$ with $\mu\in H(\la)$ can be visulized by putting
a cross ($\times$) in each cell of $\la/\mu$.

\begin{ex} For $\la=(10,9,8,6,1)$ and $\mu=(9,8,7,3,1)\in H(\la)$, 
their Ferrers diagrams and
the block $A_4(\la)$ are represented as follows~:
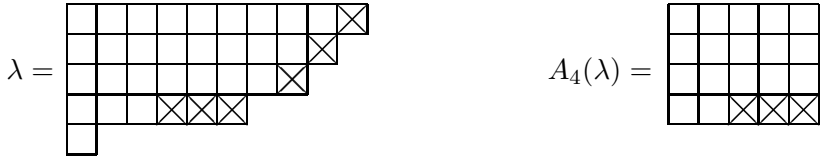
$$
\setlength{\unitlength}{0.2mm}
\begin{picture}(540,100)(60,700)
\thinlines
\put(40,750){$\la=$}\put(400,750){$A_4(\lambda)=$}
\put(100,800){\line( 0,-1){ 80}}
\put(120,720){\line( 0, 1){ 80}}
\put(140,800){\line( 0,-1){ 80}}
\put(160,720){\line( 0, 1){ 80}}
\put(180,800){\line( 0,-1){ 80}}
\put(200,740){\line( 0, 1){ 60}}
\put(220,800){\line( 0,-1){ 60}}
\put(240,760){\line( 0, 1){ 40}}
\put(260,800){\line( 0,-1){ 20}}
\put(260,780){\line(-1, 0){180}}
\put( 80,760){\line( 1, 0){160}}
\put(200,740){\line(-1, 0){120}}
\put( 80,800){\line( 1, 0){200}}
\put(280,800){\line( 0,-1){ 20}}
\put(280,780){\line(-1, 0){ 20}}
\put(260,780){\line( 0,-1){ 20}}
\put(260,760){\line(-1, 0){ 20}}
\put(260,780){\line( -1,-1){ 20}}
\put(280,800){\line( -1,-1){ 20}}
\put(240,780){\line( 1,-1){ 20}}
\put(260,800){\line( 1,-1){ 20}}
\put(220,760){\line( 1,-1){ 20}}
\put(240,760){\line( -1,-1){ 20}}
\put(200,740){\line( -1,-1){ 20}}
\put(180,740){\line( -1,-1){ 20}}
\put(160,740){\line( -1,-1){ 20}}
\put(140,740){\line( 1,-1){ 20}}
\put(160,740){\line( 1,-1){ 20}}
\put(180,740){\line( 1,-1){ 20}}
\put(240,760){\line( 0,-1){ 20}}
\put(240,740){\line(-1, 0){ 40}}
\put(200,740){\line( 0,-1){ 20}}
\put(200,720){\line(-1, 0){100}}
\put(100,720){\line( 0,-1){ 20}}
\put(100,700){\line(-1, 0){ 20}}
\put( 80,700){\line( 0, 1){100}}
\put( 80,720){\line( 1, 0){ 20}}
\put( 480,800){\line( 0,-1){ 80}}
\put( 500,720){\line( 0, 1){ 80}}
\put(520,800){\line( 0,-1){ 80}}
\put(540,720){\line( 0, 1){ 80}}
\put(480,780){\line( 1, 0){100}}
\put(480,760){\line(1, 0){100}}
\put(480,740){\line( 1, 0){100}}
\put(560,720){\line(-1, 1){ 20}}
\put(560,720){\line( 1, 1){ 20}}
\put(540,720){\line(-1, 1){ 20}}
\put(480,800){\line( 1, 0){100}}
\put(560,800){\line( 0,-1){ 80}}
\put(480,720){\line(1, 0){100}}
\put(560,720){\line( 0, 1){ 80}}
\put(560,740){\line( -1, -1){ 20}}
\put(560,740){\line( 1, -1){ 20}}
\put(580,720){\line( 0, 1){ 80}}
\put(520,720){\line( 1, 1){ 20}}
\end{picture}
$$
\end{ex}
Similarly, 
by  (\ref{iw}) and (\ref{eq:pieri}), we see that theorem~5 is equivalent to
the following:
\begin{equation}\label{eq:combi}
f_\la(a,b,c)=\sum_{(\mu,\, \nu)\in C(\la)}a^{c(\nu)}b^{|\mu/\nu|}c^{|\la/\mu|},
\end{equation}
where $C(\la)=\{(\mu, \nu)\, |\, \mu\in H(\la), \nu\in H(\mu)\}$.

We shall compute the right-hand side of (\ref{eq:combi})
using sieve method. To this end we shall first enumerate 
a larger class of patterns whose generating function is equal to
$F_{\Gamma(\la)}(a,b,c)$.

Recall that we identify a partition $\la$ with its Ferrers diagram. 
We will say that a subset $S$ of $\N^2$ is a \emph{partition diagram} if
$\{(x-k,y)|(x,y)\in S\}$ is a Ferrers diagram  for some 
integer $k\geq 0$.  Let 
$H'(\la)$ be the set of all subsets $\mu$ of $\la$ such that
$\mu\cap A_j(\lambda)$ is a partition diagram for all $j\geq 1$ and $\la/\mu$ is 
a horizontal strip.
Define
$$
B(\la)=\{(\mu,\nu)\,|\,\mu\in H(\la),\nu \in H'(\mu)\}.
$$
Note that in the above definition, the subdiagram $\nu$ of $\la$
is not necessary a partition diagram. In this regard, the set
$C(\la)$ can be described as follows:
$$
C(\la)=\{(\mu,\nu)\in B(\la)\,|\, \nu\in H(\mu)\}.
$$
Given  $\nu\in H'(\mu)$, the $j$th row of $\nu$
is called \emph{compatible} if
$$
\forall x\geq 1,\quad (x+1;j)\in\nu \Longrightarrow (x;j)\in\nu.
$$
For $p\geq 0$ let
$B_p(\la)$ be the set of $(\mu,\nu)\in
B(\la)$ such that $\nu$ 
has at least $p$ non compatible rows.
Clearly $B_0(\la)=B(\la)$ and  $ B(\la)\setminus C(\la)=B_1(\la)$, 
in other words,
a pair $(\mu,\nu)\in B(\la)$ is an element of $C(\la)$ 
iff all the rows of $\nu$ are  compatible.
By the principle of inclusion-exclusion we obtain
\begin{equation}\label{eq:inec}
\sum_{(\mu,\nu)\in C(\la)}a^{c(\nu)}b^{|\mu/\nu|}c^{|\la/\mu|}
=\sum_{p=0}^{l(\la)}(-1)^p\sum_{(\mu,\nu)\in
B_p(\la)}a^{c(\nu)}b^{|\mu/\nu|}c^{|\la/\mu|}.
\end{equation}

Each triple $(\la, \mu,\nu)$ with
 $(\mu,\nu)\in B(\la)$
 can be visualized by
putting  a circle $\circ$ (resp. cross $\times$)
in each cell of $\mu/\nu$ (resp. $\la/\mu$).

\begin{ex} 
The following diagrams represent two triples $(\la,\mu,\nu)$~:
$$
\vspace{0.5cm}
\setlength{\unitlength}{0.2mm}
\begin{picture}(680,100)(170,700)
\thinlines
\put(240,670){(a)}
\put(580,670){(b)}
\put(270,730){\circle{14}}\put(330,750){\circle{14}}
\put(370,770){\circle{14}}\put(390,790){\circle{14}}
\put(530,730){\circle{14}}\put(250,730){\circle{14}}
\put(550,730){\circle{14}}\put(570,730){\circle{14}}
\put(670,770){\circle{14}}\put(690,790){\circle{14}}
\put(220,800){\line( 1, 0){200}}
\put(420,800){\line( 0,-1){ 20}}
\put(420,780){\line(-1, 0){ 20}}\put(400,780){\line( 0,-1){ 20}}
\put(400,760){\line(-1, 0){ 20}}\put(380,760){\line( 0,-1){ 20}}
\put(380,740){\line(-1, 0){ 40}}\put(340,740){\line( 0,-1){ 20}}
\put(340,720){\line(-1, 0){100}}\put(240,720){\line( 0,-1){ 20}}
\put(240,700){\line(-1, 0){ 20}}\put(220,700){\line( 0, 1){100}}
\put(220,780){\line( 1, 0){180}}\put(220,760){\line( 1, 0){160}}
\put(340,740){\line(-1, 0){120}}\put(220,720){\line( 1, 0){ 20}}
\put(240,720){\line( 0, 1){ 80}}
\put(260,800){\line( 0,-1){ 80}}
\put(280,720){\line( 0, 1){ 80}}\put(300,800){\line( 0,-1){ 80}}
\put(320,720){\line( 0, 1){ 80}}\put(340,800){\line( 0,-1){ 60}}
\put(360,740){\line( 0, 1){ 60}}\put(380,800){\line( 0,-1){ 40}}
\put(400,780){\line( 0, 1){ 20}}\put(520,780){\line( 1, 0){180}}
\put(230,730){\circle{14}}
\put(680,760){\line(-1, 0){160}}\put(720,800){\line(-1,-1){ 20}}
\put(520,740){\line( 1, 0){120}}\put(520,720){\line( 1, 0){ 20}}
\put(540,720){\line( 0, 1){ 80}}\put(560,800){\line( 0,-1){ 80}}
\put(580,720){\line( 0, 1){ 80}}\put(600,800){\line( 0,-1){ 80}}
\put(620,720){\line( 0, 1){ 80}}\put(640,800){\line( 0,-1){ 60}}
\put(660,740){\line( 0, 1){ 60}}\put(680,800){\line( 0,-1){ 40}}
\put(700,780){\line( 0, 1){ 20}}\put(520,800){\line( 1, 0){200}}
\put(720,800){\line( 0,-1){ 20}}\put(720,780){\line(-1, 0){ 20}}
\put(700,780){\line( 0,-1){ 20}}\put(700,760){\line(-1, 0){ 20}}
\put(680,760){\line( 0,-1){ 20}}\put(680,740){\line(-1, 0){ 40}}
\put(640,740){\line( 0,-1){ 20}}\put(640,720){\line(-1, 0){100}}
\put(540,720){\line( 0,-1){ 20}}\put(540,700){\line(-1, 0){ 20}}
\put(520,700){\line( 0, 1){100}}\put(220,720){\line( 1,-1){ 20}}
\put(240,720){\line(-1,-1){ 20}}\put(280,740){\line( 1,-1){ 20}}
\put(300,740){\line(-1,-1){ 20}}\put(320,720){\line(-1, 1){ 20}}
\put(320,740){\line(-1,-1){ 20}}\put(320,740){\line( 1,-1){ 20}}
\put(340,740){\line(-1,-1){ 20}}\put(360,760){\line( 1,-1){ 20}}
\put(380,760){\line(-1,-1){ 20}}\put(380,780){\line( 1,-1){ 20}}
\put(400,780){\line(-1,-1){ 20}}\put(400,800){\line( 1,-1){ 20}}
\put(420,800){\line(-1,-1){ 20}}\put(520,720){\line( 1,-1){ 20}}
\put(540,720){\line(-1,-1){ 20}}\put(580,740){\line( 1,-1){ 20}}
\put(600,740){\line(-1,-1){ 20}}\put(620,720){\line(-1, 1){ 20}}
\put(620,740){\line(-1,-1){ 20}}\put(620,740){\line( 1,-1){ 20}}
\put(640,740){\line(-1,-1){ 20}}\put(660,740){\line( 1, 1){ 20}}
\put(660,760){\line( 1,-1){ 20}}\put(680,780){\line( 1,-1){ 20}}
\put(700,780){\line(-1,-1){ 20}}\put(700,800){\line( 1,-1){ 20}}
\end{picture}
$$
Clearly $\la=(10,9,8,6,1)$ and $\mu=(9,8,7,3)$.
In (a), the pair $(\mu,\nu)$ is in $B_1(\la)$
because the third row of $\nu$ is not compatible, so 
$\nu\in H'(\mu)\setminus 
H(\mu)$ and $\nu$ is not a partition.
In (b),
the pair $(\mu,\nu)$ is in $C(\la)$ because all the rows of $\nu$ are 
compatible, so $\nu=(8,7,7)$ is a partition in $H(\mu)$.
\end{ex}

\begin{lem} We have
$$
\sum_{(\mu,\nu)\in B(\la)}
a^{c(\nu)}b^{|\mu/\nu|}c^{|\la/\mu|}=F_{\Gamma(\la)}(a,b,c).
$$
\end{lem}
\begin{proof}
As in the proof of
(\ref{iw'}), we divide the diagram $\la$ into rectangular
blocks $A_j(\la)$, $j\geq 1$ , and
compute the weight within each block $A_j(\la)$.
Clearly choosing a pair $(\mu,\nu)$ in $B(\la)$ is equivalent to, 
for each $j\geq 1$,
first choose the $p$  left-most (resp. the rest $q=c_j-p$)
 columns of length $j$ (resp. $j-1$)
for $\mu$ in $A_j(\la)$, and then choose $s$ (resp. $p-s$)
left-most columns of length $j$ (resp. $j-1$) for $\nu$ among the $p$ 
columns of $\mu$, also choose $r$ (resp. the rest 
$q-r$) left-most columns  of length $j-1$ ( resp. $j-2$) for $\nu$.
Thus the corresponding weight is $h_{c_1}(a,b,c)$ if  $j=1$
and, for each $j\geq 2$,
$$
\left\lbrace\begin{array}{ll}\displaystyle
\sum_{p+q=c_j}c^q\left(\sum_{s=0}^p(ab)^s\right)
\left(\sum_{r=0}^{q}a^rb^{q-r}\right)
=P_{c_j}(a,b,c)&\;\textrm{ if  $j$ even};\\
&\\\displaystyle
\sum_{p+q=c_j}c^q\left(\sum_{s=0}^pb^s\right)
\left(\sum_{r=0}^q(ab)^r\right)=Q_{c_j}(a,b,c) &\; \textrm{if $j$ odd}.
	      \end{array}\right.
$$
Multiplying up over all $j\geq 1$ we get the desired formula.
\end{proof}

\begin{ex} Consider the (a) case of 
the previous example.  The subdiagrams corresponding to the block $A_4(\la)$ are
the following~:

\vspace{0.5cm}
\setlength{\unitlength}{0.2mm}
\begin{picture}(440,160)(40,640)
\thinlines
\put(40,715){$A_4(\lambda)=$}
\put(260,715){$\longrightarrow\hspace{1cm}\left\{\raisebox{-1.2cm}{}\right.$}
\put(370,760){$\mu\cap A_4(\la)=$}
\put(370,665){$\nu\cap A_4(\la)=$}
\put(130,690){\circle{14}}\put(210,710){\circle{14}}
\put(480,800){\line( 0,-1){ 80}}\put(480,720){\line( 1, 0){ 40}}
\put(520,720){\line( 0, 1){ 20}}\put(520,740){\line( 1, 0){ 60}}
\put(580,740){\line( 0, 1){ 60}}\put(580,800){\line(-1, 0){100}}
\put(480,700){\line( 1, 0){100}}\put(580,700){\line( 0,-1){ 40}}
\put(580,660){\line(-1, 0){ 20}}\put(560,660){\line( 0,-1){ 20}}
\put(560,640){\line(-1, 0){ 80}}\put(480,640){\line( 0, 1){ 60}}
\put(120,760){\line( 0,-1){ 80}}\put(120,680){\line( 1, 0){100}}
\put(220,680){\line( 0, 1){ 80}}\put(220,760){\line(-1, 0){100}}
\put(140,760){\line( 0,-1){ 80}}\put(160,680){\line( 0, 1){ 80}}
\put(180,760){\line( 0,-1){ 80}}\put(200,680){\line( 0, 1){ 80}}
\put(220,740){\line(-1, 0){100}}\put(120,720){\line( 1, 0){100}}
\put(220,700){\line(-1, 0){100}}\put(500,800){\line( 0,-1){ 80}}
\put(500,700){\line( 0,-1){ 60}}\put(520,640){\line( 0, 1){ 60}}
\put(150,690){\circle{14}}\put(520,740){\line( 0, 1){ 60}}
\put(180,700){\line(-1,-1){ 20}}\put(540,800){\line( 0,-1){ 60}}
\put(540,700){\line( 0,-1){ 60}}\put(560,660){\line( 0, 1){ 40}}
\put(560,740){\line( 0, 1){ 60}}\put(580,780){\line(-1, 0){100}}
\put(480,760){\line( 1, 0){100}}\put(520,740){\line(-1, 0){ 40}}
\put(480,680){\line( 1, 0){100}}\put(560,660){\line(-1, 0){ 80}}
\put(220,700){\line(-1,-1){ 20}}\put(200,700){\line( 1,-1){ 20}}
\put(180,700){\line( 1,-1){ 20}}\put(200,700){\line(-1,-1){ 20}}
\put(160,700){\line( 1,-1){ 20}}
\end{picture}\\
Note that $c_j=5$, $p=2$, $s=0$ and $r=2$.
\end{ex}

For any  set of integers $J=\{j_1, j_2,\ldots, j_p\}$ ($p\geq 1$) let
$B_{J}(\la)$ denote the set of all the pairs
$(\mu,\nu)\in B(\la)$ such that the $j$th row 
of $\nu$ is not compatible for $j\in J$. Hence $B_J(\la)\in B_p(\la)$.
\begin{lem} There holds
$$ 
\sum_{(\mu,\nu)\in
B_J(\la)}a^{c(\nu)}b^{|\mu/\nu|}c^{|\la/\mu|}
=(abc)^p\delta_{j_1}\ldots
\delta_{j_p}F_{\Gamma(\la)}(a,b,c).
 $$
\end{lem}
\begin{proof} Recall that $\la'=(1^{c_1}2^{c_2}\ldots)$. 
Suppose there exists a pair $(\mu,\nu)$ in $B_J(\la)$, then there should be an integer
$x_j\in \N$ such that  
$(x_j+1, j)\in \nu$ and 
$ (x_j, j)\in \mu/\nu$ for any $j\in J$. 
In view of the definition of $B_J(\la)$
we must have $x_j=c_{l(\la)}+\cdots +c_{j+1}$ and $(x_j,j+1)\in \la/\mu$,
for $\la/\mu$ is a horizontal strip.
It follows that 
$c_j\geq 1$ and $c_{j+1}\geq  1$. 
Furthermore, if $j+1$ is also in $J$, 
we must have $c_{j+1}\geq 2$.
Summarizing, we have the following equivalence:
$$
B_J(\la)\neq \emptyset \Longleftrightarrow 
c_jc_{j+1}\neq 0\; \forall j\in J\;\textrm{and }\;
c_{j+1}\geq 2 \;\textrm{if}\;  j, j+1\in J.
$$
It is easy to see that the last condition is
equivalent to $\delta_{j_1}\cdots \delta_{j_p}\Gamma(\la)\in \N^\infty$
or $\delta_{j_1}\cdots \delta_{j_p}F_{\Gamma(\la)}(a,b,c)\neq 0$.

In what follows we  shall assume that $B_J(\la)\neq \emptyset$.
Thus  we can define a unique partition 
$\delta_J(\la)$ such that
$\Gamma(\delta_J(\la))=\delta_{j_1}\ldots \delta_{j_p}\Gamma(\la)$.
Graphically, the diagram $\delta_J(\la)$ can be  obtained
 by deleting,  successively 
 for $j\in J$, the $x_j$th and $(x_j+1)$th
columns and shift all the cells on the right of $x_j$th column
 of $\la$  to left by two units. 
For $(\mu,\nu)\in B_J(\la)$, if we 
apply the same graphical operation to 
the $\mu$ and $\nu$, we get a pair $(\delta_J(\mu), \delta_J(\nu))\in 
B(\delta_J(\la))$.\\
For example, in the previous example, if 
$J=\{3,4\}$, then
$\delta_J(\la)=(6,5,4,3)$. The corresponding 
 triples $(\la,\mu,\nu)$ with $(\mu,\nu)\in B(\la)$ and 
 $(\delta_J(\la),\delta_J(\mu), \delta_J(\nu))$ with 
$(\delta_J(\mu), \delta_J(\nu))\in
B(\delta_J(\la))$ are illustrated as follows~:

\vspace{1cm}
\setlength{\unitlength}{0.2mm}
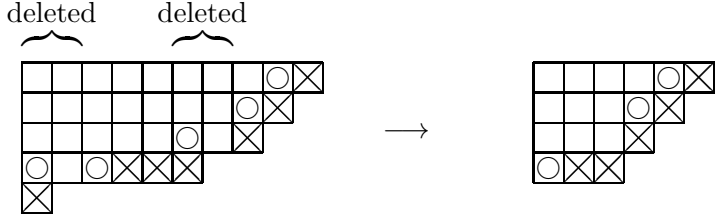
\begin{picture}(300,100)(0,700)
\thinlines
\put(280,750){$\longrightarrow$}
\put(30,810){$\overbrace{\hspace{0.8cm}}^{\hbox{deleted}}$}
\put(130,810){$\overbrace{\hspace{0.8cm}}^{\hbox{deleted}}$}
\put(150,750){\circle{14}}
\put(190,770){\circle{14}}
\put(210,790){\circle{14}}
\put( 90,730){\circle{14}}
\put(390,730){\circle{14}}
\put(450,770){\circle{14}}
\put(470,790){\circle{14}}
\put( 40,800){\line( 1, 0){200}}
\put(240,800){\line( 0,-1){ 20}}
\put(240,780){\line(-1, 0){ 20}}
\put(220,780){\line( 0, 1){  0}}
\put(220,780){\line( 0,-1){ 20}}
\put(220,760){\line(-1, 0){ 20}}
\put(200,760){\line( 0,-1){ 20}}
\put(200,740){\line(-1, 0){ 40}}
\put(160,740){\line( 0,-1){ 20}}
\put(160,720){\line(-1, 0){100}}
\put( 60,720){\line( 0,-1){ 20}}
\put( 60,700){\line(-1, 0){ 20}}
\put( 40,700){\line( 0, 1){100}}
\put( 60,800){\line( 0,-1){ 80}}
\put( 80,800){\line( 0,-1){ 80}}
\put(100,800){\line( 0,-1){ 80}}
\put(120,800){\line( 0,-1){ 80}}
\put(140,800){\line( 0,-1){ 80}}
\put(160,800){\line( 0,-1){ 60}}
\put(180,800){\line( 0,-1){ 60}}
\put(200,800){\line( 0,-1){ 40}}
\put(220,800){\line( 0,-1){ 20}}
\put( 40,780){\line( 1, 0){180}}
\put( 40,760){\line( 1, 0){160}}
\put( 40,720){\line( 1, 0){ 20}}
\put( 40,720){\line( 1,-1){ 20}}
\put( 60,720){\line(-1,-1){ 20}}
\put(100,740){\line( 1,-1){ 20}}
\put(120,740){\line(-1,-1){ 20}}
\put(120,740){\line( 1,-1){ 20}}
\put( 40,740){\line( 1, 0){120}}
\put( 50,730){\circle{14}}
\put(140,740){\line(-1,-1){ 20}}
\put(500,800){\line(-1,-1){ 20}}
\put(140,740){\line( 1,-1){ 20}}
\put(160,740){\line(-1,-1){ 20}}
\put(180,760){\line( 1,-1){ 20}}
\put(200,760){\line(-1,-1){ 20}}
\put(200,780){\line( 1,-1){ 20}}
\put(220,780){\line(-1,-1){ 20}}
\put(220,800){\line( 1,-1){ 20}}
\put(240,800){\line(-1,-1){ 20}}
\put(380,800){\line( 0,-1){ 80}}
\put(380,720){\line( 1, 0){ 60}}
\put(440,720){\line( 0, 1){ 20}}
\put(440,740){\line( 1, 0){ 20}}
\put(460,740){\line( 0, 1){ 20}}
\put(460,760){\line( 1, 0){ 20}}
\put(480,760){\line( 0, 1){ 20}}
\put(480,780){\line( 1, 0){ 20}}
\put(500,780){\line( 0, 1){ 20}}
\put(500,800){\line(-1, 0){120}}
\put(380,780){\line( 1, 0){100}}
\put(460,760){\line(-1, 0){ 80}}
\put(380,740){\line( 1, 0){ 60}}
\put(440,740){\line( 0, 1){ 60}}
\put(460,800){\line( 0,-1){ 40}}
\put(480,780){\line( 0, 1){ 20}}
\put(420,800){\line( 0,-1){ 80}}
\put(400,720){\line( 0, 1){ 80}}
\put(420,740){\line( 1,-1){ 20}}
\put(440,740){\line(-1,-1){ 20}}
\put(400,740){\line( 1,-1){ 20}}
\put(420,740){\line(-1,-1){ 20}}
\put(440,760){\line( 1,-1){ 20}}
\put(460,760){\line(-1,-1){ 20}}
\put(460,780){\line( 1,-1){ 20}}
\put(480,780){\line(-1,-1){ 20}}
\put(480,800){\line( 1,-1){ 20}}
\end{picture}\\ 
Since the weight
 corresponding to the deleted $x_j$th and $x_{j+1}$th
columns of $\la$, $\mu$ and $\nu$ is $abc$ for each $j\in J$,  we have
$$
a^{c(\nu)}b^{|\la/\nu|}c^{|\mu/\la|}=(abc)^p
a^{c(\delta_J(\nu))}b^{|\delta_J(\la)/\delta_J(\nu|)}
c^{|\delta_J(\mu)/\delta_J(\la|)}.
$$
Therefore
$$
\sum_{(\la,\nu)\in
B_J(\la)}a^{c(\nu)}b^{|\la/\nu|}c^{|\mu/\la|}
=(abc)^p\sum_{(\la,\nu)\in
B(\delta_J(\la))}a^{c(\nu)}b^{|\la/\nu|}c^{|\mu/\la|}.
$$
The lemma follows then immediately from lemma~1.
\end{proof}

It follows from lemma~2 that for $p\geq 1$
$$
\sum_{(\mu,\nu)\in
B_p(\la)}a^{c(\nu)}b^{|\mu/\nu|}c^{|\la/\mu|}
=(abc)^p\sum_{1\leq j_1<\cdots <j_p\leq l(\la)}
\delta_{j_1}\ldots \delta_{j_p}F_{\Gamma(\la)}(a,b,c).
$$
Combining  with (\ref{eq:inec}) and lemma~1 we derive immediately theorem~1.

\bigskip 
\begin{rem}
Similarly,  using another identity of Littlewood~\cite{Ma}:
\begin{equation}\label{lit'}
\sum_{\la}a^{r(\la)}s_\la(X)=\prod_{i}\frac{1+ax_i}{1-x_i^2}\prod_{j<k}
(1-x_jx_k)^{-1},
\end{equation}
where $r(\la)$ is the number of rows of odd length of $\la$, 
we obtain:
$$ \sum_{\la}f_{\la'}(a,b)\,s_\la(X)
=\prod_{i}\frac{(1+ax_i)(1+bx_i)}{1-x_i^2}\prod_{j<k}(1-x_jx_k)^{-1}
$$ and $$ \sum_{\la}f_{\la'}(a,b,c)\,s_\la(X)
=\prod_{i}\frac{(1+ax_i)(1+bx_i)
(1+cx_i)}{1-x_i^2}\prod_{j<k}(1-x_jx_k)^{-1}.
$$ 
Note  that $c_j(\la')=m_j(\la)$ is the
multiplicity of $j$ in $\la$.
\end{rem}

\section{Proof of theorem~5}
Consider the generating function $$
S(u)=\sum_{\la_0,\la}f_\la(a,b)\,s_\la(X)\, u^{\la_0} $$
where the sum is over all
$\la_0\geq \la_1\geq \cdots \geq \la_n\geq 0$, and 
$\la=(\la_1,\ldots, \la_n)$. 
 Suppose
$\la$ is of form $\mu_1^{r_1},\, \mu_2^{r_2},\, \ldots \mu_k^{r_k}$,
 where
$\mu_1>\mu_2>\cdots >\mu_k\geq 0$ and the $r_i$ are positive
integers whose sum is $n$.
Let $S_n^\la=S_{r_1}\times \cdots
S_{r_k}$ be the group of permutations leaving $\la$ invariant. Then
\begin{eqnarray*}
s_\la(X)&=&\sum_{w\in
S_n}w\left(x_1^{\la_1}\ldots x_n^{\la_n}
\prod_{i<j}\frac{x_i}{ x_i-x_j}\right)\\
&=& \sum_{w\in
S_n/S_n^\la}w\left(x_1^{\la_1}\ldots x_n^{\la_n}
\prod_{\la_i>\la_j}\frac{x_i}{ x_i-x_j}\right),
\end{eqnarray*}
where the permutation $w$ acts on the indices of the
indeterminates. Each 
$w\in S_n/S_n^\la$ corresponds to  a surjective mapping $f:
X\longrightarrow \{1,2,\ldots, k\}$ such that $|f^{-1}(i)|=r_i$.
For any subset $Y$ of $X$, let $p(Y)$ denote  the product of the 
elements of $Y$. (In particular, $p(\emptyset)=1$.)
We can rewrite Schur functions as follows:
$$
s_\la(X)=\sum_{f}
p(f^{-1}(1))^{\mu_1}\cdots p(f^{-1}(k))^{\mu_k}
\prod_{f(x_i)<f(x_j)}{x_i\over x_i-x_j}. 
$$
summed over all  surjective mappings $f:
X\longrightarrow \{1,2,\ldots, k\}$ such that $|f^{-1}(i)|=r_i$.
Furthermore,
each such $f$ determines a \emph{filtration}  of $X$:
$$
{\cal F}: \quad \emptyset=F_0\subsetneq F_1\subsetneq 
\cdots \subsetneq F_k=X, 
$$
according to the rule $x_i\in F_l\Longleftrightarrow f(x_i)\leq l$
for $1\leq l\leq k$.
Conversely, such a filtration
${\cal F}=(F_0,\, F_1, \ldots, F_k)$ determines a surjection 
$f: X\longrightarrow \{1,2,\ldots, k\}$ uniquely.
Thus we can write:
\begin{equation}\label{filtre}
s_\la(X)=\sum_{\cal F}\pi_{\cal F}\prod_{1\leq i\leq
k}p(F_i\setminus F_{i-1})^{\mu_i},
\end{equation}
summed  over all the filtrations $\cal F$  such that
 $|F_i|=r_1+r_2+\cdots +r_i$ for $1\leq i\leq k$, and
$$ 
\pi_{\cal F}=\prod_{f(x_i)<f(x_j)}{x_i\over x_i-x_j}, 
$$ 
where $f$ is the function defined by $\cal F$.

Now let
$\nu_i=\mu_i-\mu_{i+1}$ if $1\leq i\leq k-1$ and $\nu_k=\mu_k$,
thus $\nu_i>0$ if $i<k$ and $\nu_k\geq 0$. 
Since the  lengths of columns of $\la$ are 
 $|F_j|=r_1+\cdots r_j$  
with multiplicities $\nu_j$ for $1\leq j\leq k$, we have 
\begin{equation}\label{eq:mac}
f_\la(a,b)=\prod_{ |F_j| \,
\hbox{\scriptsize odd}}{a^{\nu_j+1}-b^{\nu_j+1}\over a-b}\,
\prod_{|F_j| \, \hbox{\scriptsize even}}{1-(ab)^{\nu_j+1}\over
1-ab}. 
\end{equation}
Furthermore,
let  $\mu_0=\la_0$ and $\nu_0=\mu_0-\mu_1$ in the definition of
$S(u)$, so that $\nu_0\geq 0$ and  $\mu_0=\nu_0+\nu_1+\cdots +\nu_k$.
 It follows from 
(\ref{filtre}) and (\ref{eq:mac}) that~:
\begin{eqnarray}
S(u) &=& \sum_{\cal F}\pi_{\cal F}\sum_{\nu}u^{\nu_0}
\prod_{|F_j|\, \hbox{\scriptsize odd}}{a^{v_j+1}-b^{v_j+1}\over
a-b} u^{v_j}p(F_j)^{v_j}\,\nonumber \\
&&\hskip 2 cm \times \prod_{|F_j|\,
\hbox{\scriptsize even}} {1-(ab)^{v_j+1}\over
1-ab}u^{v_j}p(F_j)^{v_j},\label{eq:Mac}
\end{eqnarray}
where the outer sum is over all filtrations $\cal F$ of $X$ 
and the inner sum is over all integers $\nu_0,\nu_1,\ldots, \nu_k$ such that 
$\nu_0\geq 0$, $\nu_k\geq 0$ and 
$\nu_i\geq 0$ for $1\leq i\leq k-1$.
For any filtration $\cal F$ of $X$ set
\begin{eqnarray*}
{\cal A}_{\cal F}(X, u)&=&
\prod_{|F_j|\, \hbox{\scriptsize odd}} \left[
{a(a-b)^{-1}\over 1-ap(F_j)u}-{b(a-b)^{-1}\over 1-bp(F_j)u}-
\chi(F_j\neq X)\right]\\ 
&&\times 
\prod_{|F_j|\, \hbox{\scriptsize
even}}\left[{(1-ab)^{-1}\over 1-p(F_j)u}- {ab(1-ab)^{-1}\over
1-abp(F_j)u}-\chi(F_j\neq X)\right],
\end{eqnarray*}
where $\chi(A)=1$ if $A$ is true, and
$\chi(A)=0$ if $A$ is false. Then the inner sum of (\ref{eq:Mac}) is 
$$
(1-u)^{-1}{\cal A}_{\cal F}(X, u),
$$
therefore
$$
S(u)=(1-u)^{-1}\sum_{\cal F}\pi_{\cal F}{\cal A}_{\cal F}(X, u),
$$
where the sum is over all the filtrations of $X$ as before.

The above formula shows that $S(u)$ is a rational function of $u$ whose
denominator is the product of the form $1-ap(Y)u$, $1-bp(Y)u$ or 
$1-abp(Y)u$, where
$Y\subseteq X$.
Therefore  we have the
following result.
\begin{lem} The generating function $S(u)$ is of the form:
\begin{eqnarray*}
S(u)&=&{c(\emptyset)\over 1-u}+\sum_{{Y\subseteq X}\atop |Y| \,
\hbox{\scriptsize odd}} \left({a(Y)\over 1-ap(Y)u}-{b(Y)\over
1-bp(Y)u}\right)\\ 
&&+\sum_{{Y\subseteq X}\atop |Y| \,
\hbox{\scriptsize even}>0} \left({c(Y)\over 1-p(Y)u}-{d(Y)\over
1-abp(Y)u}\right).
\end{eqnarray*}
\end{lem}

It remains to compute the residues.
Let us  start with $c(\emptyset)$. Writing $\la_0=\la_1+k$ with
$k\geq 0$, we see that
\begin{eqnarray*}
S(u)&=&\sum_{k\geq 0}u^k\sum_\la f_\la(a,b)\,
s_\la(X)u^{\la_1}\\ 
&=&(1-u)^{-1}\sum_\la f_\la(a,b)\, s_\la(X)u^{\la_1},
\end{eqnarray*}
it follows from (2) that
$$
c(\emptyset)=\left(S(u)(1-u)\right)|_{u=1}=\Phi(X;a,b).
$$
For computations of the other  residues, we introduce some more notations. For 
any $Y\subseteq X$,
let $Y'=X\setminus Y$ and $-Y=\{x_i^{-1}:x_i\in Y\}$. 
For any subset $Z$ of 
$X$ or $-X$ let
$$
\alpha(Z,u)=\left\lbrace\begin{array}{ll}
{(1-ap(Z)u)(1-bp(Z)u)} \quad&\textrm{ if $|Z|$ odd};\\
{(1-p(Z)u)(1-abp(Z)u)}
\quad&\textrm{if $|Z|$ even.}
\end{array}\right.
$$
As the computations of other residues are similar, 
we just give the 
details for $c(Y)$. Let $Y\subseteq X$ such that  $|Y|$ is even. 
Then we have 
\begin{equation}\label{eq:residu}
c(Y)=\left[(1-u)^{-1}\sum_{\cal F}\pi_{\cal F}{\cal
A}_{\cal F}(X;u)(1-p(Y)u)\right]_{u=p(-Y)}.
\end{equation}
If $Y\notin\cal F$, the corresponding summand is  equal to  0. Thus we need
only to consider the following filtrations ${\cal F}$:
$$  \emptyset=F_0\subsetneq
\cdots \subsetneq F_t=Y\subsetneq \cdots \subsetneq  F_k=X\qquad 1\leq t\leq k.
$$ 
We may then split $\cal F$ into two filtrations
${\cal F}_1$ and ${\cal F}_2$, of $-Y$ and $Y'=X\setminus
Y$ respectively, as follows :
\begin{eqnarray*}
{\cal F}_1&:& \emptyset \subsetneq  -(Y\setminus F_{t-1})\subsetneq 
\cdots \subsetneq -(Y\setminus F_1)\subsetneq  -Y,\\ 
{\cal F}_2&:&
\emptyset \subsetneq F_{t+1}\setminus Y\subsetneq   \cdots \subsetneq 
F_{k-1}\setminus Y\subsetneq  Y'.
\end{eqnarray*}
Then, writing  $v=p(Y)u$, we have  
\begin{eqnarray*}
(1-u)^{-1}{\cal A}_{\cal
F}(X;u)(1-p(Y)u)&=&{(1-p(-Y)v)^{-1}{\cal A}_{{\cal F}_1}(-Y;v)
{\cal A}_{{\cal F}_2}(Y';v)}\\
&&\times \alpha(-Y,v)\left[(1-ab)^{-1}-\beta(v)(1-v)\right],
\end{eqnarray*}
where $\beta(v)=ab/(1-abv)(1-ab)-\chi(Y\neq X)$, and 
$$ \pi_{\cal
F}(X)=\pi_{{\cal F}_1}(-Y)\pi_{{\cal F}_2}(Y')\prod_{x_i\in Y, x_j\in Y'}
(1-x_i^{-1}x_j)^{-1}, 
$$ 
As  $u=p(-Y)$ is equivalent to $v=1$, it follows from (\ref{eq:residu}) that
\begin{eqnarray*}
c(Y)&=&(1-ab)^{-1}(1-p(-Y))^{-1}\alpha(-Y,1)\prod_{x_i\in Y, x_j\in Y'}
(1-x_i^{-1}x_j)^{-1}\\ 
&&\times \left[\sum_{{\cal F}_1}
\pi_{{\cal F}_1}(-Y){\cal A}_{{\cal F}_1}(-Y;v)\right]_{v=1}
\times \left[\sum_{{\cal F}_2}
\pi_{{\cal F}_2}(Y'){\cal A}_{{\cal F}_2}(Y';v)\right]_{v=1}.
\end{eqnarray*}
Using the result of $c(\emptyset)$, which can be written: 
$$
\Phi(X;a,b)=\sum_{{\cal F}}
\left(\pi_{{\cal F}}(X)
{\cal A}_{{\cal F}}(X;u)\right)_{u=1},
$$ 
we obtain:
$$
c(Y)={\alpha(-Y,1)\Phi(-Y;a,b)\Phi(Y';a,b)\over (1-ab)(1-p(-Y))}
\prod_{x_i\in Y, x_j\in Y'}
(1-x_i^{-1}x_j)^{-1}.
$$
Each subset $Y$ of $X$ can be encoded by  a sequence $\xi\in \{\pm 1\}^n$
according to the rule~: 
$\xi_i=1$ if $x_i\notin Y$ and $\xi_i=-1$ if $x_i\in Y$.
Hence
$$
c(Y)={\Phi(x_1^{\xi_1}, \ldots, x_n^{\xi_n};a,b)\over
(1-ab)(1-p(-Y))}\alpha(-Y,1),
$$
Note also that 
$$ 
p(Y)=\prod_i x_i^{(1-\xi_i)/2},\qquad
p(-Y)=\prod_i x_i^{(\xi_i-1)/2}.
$$
In the same way, we find for any even size subset $Y\subseteq  X$ that
$$
d(Y)={ab\Phi(x_1^{\xi_1}, \ldots, x_n^{\xi_n};a,b)\over
(1-ab)(1-(ab)^{-1}p(-Y))}\alpha(-Y,1),
$$
and for any odd size subset $Y\subseteq  X$ that
\begin{eqnarray*}
a(Y)&=&{a\Phi(x_1^{\xi_1}, \ldots, x_n^{\xi_n};a,b)\over
(a-b)(1-a^{-1}p(-Y))}\alpha(-Y,1),\\ 
b(Y)&=&{b\Phi(x_1^{\xi_1},
\ldots, x_n^{\xi_n};a,b)\over (a-b)(1-b^{-1}p(-Y))}\alpha(-Y,1).
\end{eqnarray*}

By virtue of lemma~1, extracting the coefficient of $u^m$ in $S(u)$ yields
\begin{eqnarray*}
\sum_{\la\subseteq (m^n)}f_\la(a,b)s_\la(X)=  \Phi(X;a,b)&+&
\sum_{{Y\subseteq X}\atop |Y| \, \hbox{\scriptsize
odd}}\left[a(Y)a^m-b(Y)b^m\right]p(Y)^m\\ &+&\sum_{{Y\subseteq
X}\atop |Y| \, \hbox{\scriptsize even}>0}
\left[c(Y)-d(Y)(ab)^m\right]p(Y)^m.
\end{eqnarray*}
Finally,
substituting  the values of $a(Y), b(Y), c(Y)$ and $d(Y)$ in the above formula
we obtain 
theorem~5.

\section{Three special cases}
First we note that $f_\la(1,0)=1$,
$$
f_\la(1,-1)=\left\lbrace\begin{array}{ll}
0&\quad \textrm{if any $c_j$ is odd},\\
1&\quad \textrm{otherwise};
			\end{array}\right.
$$
and
$$
f_\la(0,0)=\left\lbrace\begin{array}{ll}
0&\quad \textrm{if any $c_j$ is positive for any  odd $j$},\\
1&\quad \textrm{otherwise}.
			\end{array}\right.
$$
On the other hand, we have $$\beta(\xi, 1,0)=1,$$  
$$
\beta(\xi, 1,-1)=\left\lbrace\begin{array}{ll}
1&\quad \textrm{if $m$ even},\\
\prod_ix_i^{(\xi_i-1)/2}&\quad \textrm{if $m$ odd};
			\end{array}\right.
$$
and
$$
\beta(\xi, 0,0)=\left\lbrace\begin{array}{ll}
0&\quad \textrm{if any $|\xi|_{-1}$ is odd},\\
1&\quad \textrm{otherwise}.
			\end{array}\right.
$$ 
So we derive immediately from theorem~5 the following result.
\begin{coro} 
The sums of Schur functions of shape in a given rectangle are:
\begin{eqnarray} 
\sum_{\la\subseteq (m)^n}s_\la(X)&=& \sum_{\xi\in \{\pm1\}^n}
\Phi(X^\xi;1,0)\prod_ix_i^{m(1-\xi_i)/2}, \label{a}\\ 
\sum_{{\la\subseteq (2m)^n\atop \la\; \textrm{even}}}s_\la(X)&=& \sum_{\xi\in 
\{\pm1\}^n}
\Phi(X^\xi;1,-1)\prod_ix_i^{m(1-\xi_i)}, \label{b}\\ 
\sum_{{\la\subseteq (m)^n \atop \la'\; \textrm{even}}}s_\la(X)&=& 
\sum_{{\xi\in \{\pm1\}^n\atop |\xi|_{-1}\, \textrm{even}}}
\Phi(X^\xi;0,0)\prod_ix_i^{m(1-\xi_i)/2}, \label{c}
\end{eqnarray}
where $n$ is even in the last identity.
\end{coro}

To see that the above corollary is equivalent to
theorems 2, 3 and 4, we need only to appeal to
Vandermonde's determinantal formula ~:
\begin{equation}\label{d1}
\sum_{\sigma\in 
S_n}\epsilon(\sigma)\prod_{i=1}^nx_{\sigma(i)}^i=\det(x_j^{i-1})=\prod_{1\leq 
i<j\leq n}(x_i-x_j).
\end{equation}
Notice that for $\xi\in \{\pm 1\}^n$ and $1\leq i<j\leq n$, 
$$
(x_i^{\xi_i}-x_j^{\xi_j})(1-x_i^{\xi_i}x_j^{\xi_j})
=(x_i-x_j)(1-x_ix_j)x_i^{\xi_i-1}x_j^{\xi_j-1},
$$
therefore
\begin{equation}\label{d2} 
\prod_{i<j}(x_i^{\xi_i}-x_j^{\xi_j})(1-x_i^{\xi_i}x_j^{\xi_j})\\
=\prod_{i<j}(x_i-x_j)(1-x_ix_j)\prod_ix_i^{(n-1)(\xi_i-1)}.
\end{equation}

\textbf{The $(a,b)=(1,0)$ case}~: Set 
$$
\Delta_B=\frac{\prod_{i<j}(x_i-x_j)}{ \Phi(X; 
1,0)}=\prod_i(1-x_i)\prod_{i<j}(x_i-x_j)(x_ix_j-1).
$$
Using (\ref{d1}) and (\ref{d2}), we can write
$$
\Phi(X^\xi;1,0)=\frac{(-1)^{|\xi|_{-1}}}{\Delta_B}
\prod_ix_i^{(1-\xi_i)(n-1/2)}\sum_{\sigma\in S_n}
\epsilon(\sigma)\prod_ix_{\sigma(i)}^{\xi_{\sigma(i)}(i-1)}.
$$
So the right side of (\ref{a}) is
\begin{eqnarray*}
&&\frac{1}{\Delta_B}\sum_{\sigma\in
S_n}\epsilon(\sigma)\sum_{\xi\in \{\pm
1\}^n}(-1)^{|\xi|_{-1}}\prod_ix_{\sigma(i)}^{(m+2n-1)(1-\xi_{\sigma(i)})/2+\xi_{
\sigma(i)}(i-1)}\\
&&=\frac{1}{\Delta_B}\sum_{\sigma\in
S_n}\epsilon(\sigma)\sum_{\xi\in \{\pm
1\}^n}\prod_{\xi_{\sigma(i)}=1}x_{\sigma(i)}^{i-1}
\prod_{\xi_{\sigma(i)}=-1}\left(-x_{\sigma(i)}^{m+2n-i}\right)\\
&&=\frac{1}{\Delta_B}\det\left(x_i^{j-1}-x_i^{m+2n-j}\right).
\end{eqnarray*}
Hence theorem~2 is equivalent to (\ref{a}).

\textbf{The $(a,b)=(1,-1)$ case}~: Set 
$$
\Delta_C=\frac{\prod_{i<j}(x_i-x_j)}{\Phi(X; 1,-1)}
=\prod_i(1-x_i^2)\prod_{i<j}(x_i-x_j)(1-x_ix_j).
$$
By (\ref{d1}) and (\ref{d2}), we have also
$$\Phi(X^\xi;1,-1)=\frac{(-1)^{|\xi|_{-1}}}{\Delta_C}
\prod_ix_i^{n(1-\xi_i)}\sum_{\sigma\in S_n}
\epsilon(\sigma)\prod_ix_{\sigma(i)}^{\xi_{\sigma(i)}(i-1)},$$
and the right hand side of (\ref{b}) is
\begin{eqnarray*}
&&\frac{1}{\Delta_C}\sum_{\sigma\in
S_n}\epsilon(\sigma)\sum_{\xi\in \{\pm
1\}^n}(-1)^{|\xi|_{-1}}\prod_ix_{\sigma(i)}^{(n+m)
(1-\xi_{\sigma(i)})+(i-1)\xi_{\sigma(i)}}\\
&&=\frac{1}{\Delta_C}\sum_{\sigma\in S_n}\epsilon(\sigma)\sum_{\xi\in \{\pm
1\}^n}\prod_{\xi_{\sigma(i)}=1}x_{\sigma(i)}^{i-1}\prod_{\xi_{\sigma(i)}=-1}
\left(-x_{\sigma(i)}^{2n+2m-i+1}\right)\\
&&=\frac{1}{\Delta_C}\det\left(x_i^{j-1}-x_i^{2m+2n+1-j}\right).
\end{eqnarray*}
So Theorem~3 is equivalent to (\ref{b}).

\textbf{The $(a,b)=(0,0)$ case}~:
Set 
$$
\Delta_D=\frac{\prod_{i<j}(x_i-x_j)}{\Phi(X; 
0,0)}=\prod_{i<j}(x_i-x_j)(1-x_ix_j).
$$
By (\ref{d1}) and (\ref{d2}), we have also
$$\Phi(X^\xi;0,0)=\frac{1}{\Delta_D}\prod_ix_i^{(n-1)(1-\xi_i)}\sum_{\sigma\in 
S_n}\epsilon(\sigma)\prod_ix_{\sigma(i)}^{\xi_{\sigma(i)}(i-1)}$$
and the right side of (\ref{c}) is 
\begin{eqnarray*}
&&\frac{1}{\Delta_D}\sum_{\sigma\in
S_n}\epsilon(\sigma)\sum_{\xi\in \{\pm
1\}^n\atop |\xi|_{-1}\hbox{\scriptsize
even}}\prod_{\xi_{\sigma(i)}=1}
x_{\sigma(i)}^{i-1}\prod_{\xi_{\sigma(i)}=-1}x_{\sigma(i)}^{2n+m-i-1}\\
&&=\frac{1}{2\Delta_D}\left[\det\left(x_i^{j-1}-x_i^{m+2n-1-j}\right)+
\det\left(x_i^{j-1}+x_i^{m+2n-1-j}\right)\right].
\end{eqnarray*}
So theorem~4 is equivalent to (\ref{c}).

When $m=0$, as the left sides of (\ref{a}), (\ref{b}) and (\ref{c})
are equal to 1, we obtain the following result.
\begin{coro}
For any non negative integer $n$, we have  
\begin{eqnarray*}
\det\left(x_i^{j-1}-x_i^{2n-j}\right)&=&\prod_i(1-x_i)\prod_{i<j}(x_i-x_j)(1-x_i
x_j),\\
\det\left(x_i^{j-1}-x_i^{2n-j+1}\right)&=&\prod_i(1-x_i^2)\prod_{i<j}(x_i-x_j)(1
-x_ix_j),\\
\det\left(x_i^{j-1}+x_i^{2n-1-j}\right)&=&2\prod_{i<j}(x_i-x_j)(1-x_ix_j).
\end{eqnarray*}
\end{coro}
These are actually
Weyl's denominator formulas  for root
systems of type $B_n$, $C_n$ and $D_n$ (\cite{Br3}, p. 68-69) respectively.



\begin{thebibliography}{9}
\bibitem{Br1}\textsc{Bressoud} (D.),
\emph{Elementary proof of MacMahon's conjecture}, 
J. Alg. Combin. {\bf 7}, 253-257, 1998.

\bibitem{Br2}\textsc{Bressoud} (D.),
\emph{Elementary proofs of identities for Schur functions and
plane partitions}, The Ramanujan J., {\bf 4}, 69-80, 2000.

\bibitem{Br3}\textsc{Bressoud} (D.),
\emph{Proofs and Confirmations, 
The Story of the Alternating Sign Matrix Conjecture}, Cambridge Univ. Press, 
1999.

\bibitem{De}\textsc{D\'esarm\'enien} (J.),
\emph{La d\'emonstration des identit\'es de Gordon et MacMahon 
et de deux identit\'es nouvelles}, Strasbourg, Publ.I.R.M.A., 
Actes du 15\`eme S\'eminaire Lotharingien de Combinatoire, 
340/S-15, 39-49, 1987.

\bibitem{IW} \textsc{Ishikawa} (M.) and \textsc{Wakayama} (M.),
 \emph{Applications of Minor-Summation Formula II. 
Pfaffians and Schur Polynomials}, J. Combin. Th., Ser. A~88 (1999), 
136-157. 

\bibitem{Ma} \textsc{Macdonald} (I.G.), 
\emph{Symmetric functions and Hall
polynomials}, Clarendon Press, second edition, Oxford, 1995.

\bibitem{Ok} \textsc{Okada} (S.), \emph{Application of minor
summation formulas to rectangular-shaped representations of
classical groups}, J. Algebra {\bf 205}, 337-367, 1998.

\bibitem{St} \textsc{Stembridge} (J. R.), \emph{Hall-Littlewood
functions, plane partitions, and the Rogers-Ramanujan identities},
Trans. Amer. Math. Soc., {\bf 319}, no.2, 469-498, 1990.

\end{thebibliography}
\end{document}